\documentclass[10pt,reqno]{amsart}
\usepackage{amsfonts}
\usepackage{amsmath}
\usepackage{CJK}
\usepackage{amssymb}
\usepackage{latexsym,bm}
\usepackage[RGB]{xcolor}
\usepackage{mathrsfs}
\usepackage{amsthm}
\numberwithin{equation}{section}
\newtheorem{theorem}{Theorem}[section]

\newtheorem{lemma}[theorem]{Lemma}
\newtheorem{remark}{Remark}[section]

\newtheorem{proposition}[theorem]{Proposition}
\newcommand{\dif}{\mathrm{d}}

\begin{document}
\title[Blow-up for relaxed compressible Navier-Stokes equations]{Blow-up of solutions for relaxed compressible Navier-Stokes equations}
\author{Yuxi Hu and Reinhard Racke}
 \thanks{\noindent Yuxi Hu, Department of Mathematics, China University of Mining and Technology, Beijing, 100083 and Key Laboratory of Scientific and Engineering Computing (Ministry of Education), P.R. China, yxhu86@163.com\\
\indent Reinhard Racke, Department of Mathematics and Statistics, University of
Konstanz, 78457 Konstanz, Germany, reinhard.racke@uni-konstanz.de
}
\begin{abstract}
 We present a blow-up result for large data for relaxed compressible Navier-Stokes models avoiding the possibility of reaching the boundary of hyperbolicity. Thus a previous result is improved and further examples are given illustrating possible effects of a relaxation and contrasting the classical compressible Navier-Stokes equations without relaxation where solutions for large data exist globally. 
 {\bf AMS classification code}: 35 L 60, 35 B 44
\end{abstract}
\maketitle
\section{Introduction}
We consider the system of one-dimensional non-isentropic compressible Navier-Stokes equations,
\begin{align}\label{1.1}
\begin{cases}
 \rho_t+(\rho u)_x=0,\\
\rho  u_t+\rho u u_x+p_x=S_x,\\
E_t+(u E+pu+q-Su)_x=0.
\end{cases}
\end{align}
with $(t,x)\in \mathbb R_+\times \mathbb R$. Here, $\rho, u, p, E $ represent the fluid density, velocity,   pressure and total energy, respectively. Equations for the stress $S$ and the heat flux $q$ should be given to make the system
\eqref{1.1} closed. We shall use the following model:
\begin{align}\label{1.2}
\tau_1(\theta) (\rho q_t+\rho u\cdot q_x)+q+\kappa(\theta) \theta_x=0,
\end{align}
and
\begin{align}\label{1.3}
\tau_2 (\rho S_t+\rho u\cdot S_x)+S=\mu u_x.
\end{align}
Here $\tau_1(\theta), \tau_2>0$ are relaxation parameters, $\kappa(\theta)>0$ and $\mu>0$ denote the heat conduction and the viscosity coefficient, respectively. $\tau_2$ and $\mu$ are assumed to be constants. The constitutive equation \eqref{1.3} was proposed by Freist\"uhler \cite{Frei1, Frei2}  for the isentropic case, see also Ruggeri \cite{Rug83} and M\"uller\cite{Mu67} for a similar model in the non-isentropic case.

Furthermore, we assume that the total energy is given by
\begin{equation} \label{eq1.3a}
E=\frac{1}{2} \rho u^2 +\frac{\tau_2}{2\mu} \rho S^2+\rho e(\theta, q),
\end{equation}
 and the specific internal energy $e$ and the pressure $p$ are given by
\begin{align}\label{1.4}
e(\theta)= C_v \theta+a(\theta)q^2, \qquad\qquad  p(\rho, \theta)=R \rho \theta,
\end{align}
where
$$
a(\theta)=\frac{Z(\theta)}{\theta}-\frac{1}{2} Z^\prime(\theta)
\quad \mbox{\rm with} \quad Z(\theta)=\frac{\tau_1(\theta)}{\kappa(\theta)}.$$
 Here,   $C_v>0, R>0$ denotes the heat capacity at constant volume and the gas constant, respectively. $p$ and $e$ satisfy the usual thermodynamic equation
$$
\rho^2 e_\rho=p-\theta p_\theta.
$$

The dependence of the internal energy on $q^2$ is indicated by Coleman et al. \cite{CFO}, where they rigorously prove that for heat equations with Cattaneo-type law,  the formulation \eqref{1.4}  is consistent with the second law of thermodynamics, see also \cite{CG, CHO, TA}.

 In the constitutive relation (\ref{1.3}), in its linearized form: $\tau_2 S_t+S=\mu u_x$, the positive parameter $\tau_2$ is the relaxation time describing the time lag in the response of the stress tensor to the velocity gradient, cf. also Christov and Jordan \cite{CJo}.   Pelton et al. \cite{Peetal} showed that such a"time lag'' cannot be neglected, even for simple fluids, in
the experiments of high-frequency vibration of nano-scale mechanical devices immersed in water-glycerol mixtures. It turned out
  that, cf. also \cite{ChSa015}, equation (\ref{1.3}) provides a general formalism  to characterize the fluid-structure interaction
of nano-scale mechanical devices vibrating in simple fluids. A similar relaxed constitutive relation was already proposed by Maxwell in \cite{Max1867}, in order to describe the relation of stress tensor and velocity gradient for a non-simple fluid.

We shall consider the Cauchy problem for the functions
\[
(\rho,u,\theta,S,q): \mathbb R \times [0,+\infty)\rightarrow \mathbb R_+\times \mathbb R \times \mathbb R_+\times \mathbb R\times \mathbb R
\]
with initial conditions
\begin{align}\label{1.5}
(\rho(x,0),u(x,0),\theta(x,0),S(x,0),q(x,0))=(\rho_0,u_0,\theta_0,S_0,q_0).
\end{align}

Neglecting $\rho$ in the constitutive relations \eqref{1.2}-\eqref{1.3} and assuming $\tau_1, \kappa$ to be constants, the authors and Wang\cite{HRW} established a blow-up result under the assumption that  $(\rho-1, \theta-1, q, S)\in \Omega$ with $\Omega$ a "small'' domain. Therefore,  the solutions  in \cite{HRW}  might "blow up" in the sense that one may reach the boundary of $\Omega$. The aim of this paper is to:
 \begin{itemize}
\item establish a symmetric hyperbolic system without smallness condition,

\item  find a physical  entropy which gives lower energy estimates and some dissipation,

\item show a global existence result for small data,

\item mainly prove a blow-up of classical solutions for large data.

\end{itemize}

It should be noted that the constitutive relations \eqref{1.2}-\eqref{1.3} have many merits.  For example,  as mentioned by Freist\"uhler \cite{Frei1}, they are Galilean invariant and in a conservation form which allows one to define weak solutions. Moreover, the use of constitutive relations \eqref{1.2}-\eqref{1.3} in this paper is originally coming from the idea that putting the system into a symmetric hyperbolic system, for which the pressure $p$ should not depend on $q$ and $S$. In this regard, to satisfy the thermodynamic relation $\rho^2 e_\rho=p-\theta p_\theta$, the specific internal energy $e$ should not depend on $\rho$ (for an ideal gas). Therefore, we removed the variable $\rho$ in the formulation of $e$ in our previous paper \cite{HRW}. Then, to have an entropy equation and a "good" equation for $\theta$, we just need the new constitutive relations \eqref{1.2}-\eqref{1.3} which  coincide with the model proposed by Freist\"uhler, at least for the isentropic case.

 The most interesting aspect might be that the blow-up result contrasts the situation without relaxation. i.e. for the classical compressible Navier-Stokes system corresponding to $\tau_1=\tau_2=0$, where large global solutions exist, see Kazhikhov \cite{KA}. This really nonlinear effect -- loosing the global existence for large data --, not anticipated from the linearized version, shows the possible impact a relaxation might have. For several linear systems of various type an effect is visible in loosing exponential stability in bounded domains or becoming of regularity loss type in the Cauchy problem, see the discussion in our paper \cite{HRW}.

 The method we use to prove the blow-up result is mainly motivated by  Sideris' paper \cite{S85} where he showed that any $C^1$ solutions of compressible Euler equations must blow up in finite time.
A blow-up result for a similar system has also been stated recently by Freist\"uhler \cite{Frei3} applying the general result for symmetric hyperbolic systems with sources in one space dimension by B\"arlin \cite{Bae023}. A solution remains bounded, but the solution does not remain in $C^1$, provided the data are \textit{small} enough.
 We shall show that the system \eqref{1.1}-\eqref{1.5} is a symmetric hyperbolic system which has the important property of finite propagation speed. This allows us to define some averaged quantities (different from that in \cite{S85}) and finally show a blow-up of solutions in finite time by establishing a Riccati-type inequality.
 In contrast to \cite{Frei3,Bae023}, our blow-up requires \textit{large} initial velocities; moreover, here the largeness is described explicitly. For initial data being small in higher-order Sobolev spaces ($H^2$), there exist global solutions. The method used here also extends to higher dimensions, see \cite{HR-hd}.

The paper is organized as follows. In Section 2, we derive an entropy equation for system \eqref{1.1}-\eqref{1.5} and then present the local existence theorem  in Section 3 together with some remarks on global existence for small data. In Section 4 we show the blow-up result.

Finally, we introduce some notation.  $W^{m,p}=W^{m,p}(\mathbb R),\,0\le m\le
\infty,\,1 \le p\le \infty$, denotes the usual Sobolev space with norm $\|
\cdot \|_{W^{m,p}}$, $H^m$ and $L^p$
stand for $W^{m,2}$ resp. $W^{0,p}$.
%After some calculation, one possible way to accheive the above aim is to modify the constitutive equation as
%\begin{align}\label{1.4}
%\tau_1 (\rho q_t+\rho u\cdot q_x)+\frac{q+\kappa \theta_x}{\theta}=0,
%\end{align}
%and
%\begin{align}\label{1.5}
%\tau_2 (\rho S_t+\rho u\cdot S_x)+S=\mu u_x.
%\end{align}

%Moreover, we assume that the internal energy $e$ and the pressure $p$ satisfy the following constitutive relations,
%\begin{align}
%e=C_v \theta+\frac{\tau_1}{2\kappa } q^2+\frac{\tau_2}{2\mu}S^2, \label{1.4}\\
%p=R \rho \theta  \label{1.5},
%\end{align}
%with positive constants $C_v, R$ denoting the heat capacity at constant volume and the gas constant, respectively,
%such that they satisfy the usual thermodynamic equation
%$$
%\rho^2 e_\rho=p-\theta p_\theta.
%$$

%\begin{align}
%\begin{cases}
% \rho_t+(\rho u)_x=0,\\
%\rho  u_t+\rho u u_x+R\theta \rho_x+R \rho \theta_x=S_x,\\
%\rho C_v \theta_t+(\rho u C_v-\frac{q}{\theta})\theta_x+R \rho \theta u_x+q_x=\frac{q^2}{\kappa \theta}+\frac{S^2}{\mu}\\
%\tau_1 (\rho q_t+\rho u\cdot q_x)+\frac{q+\kappa \theta_x}{\theta}=0,\\
%\tau_2 (\rho S_t+\rho u\cdot S_x)+S=\mu u_x.
%\end{cases}
%\end{align}
%
%The above system is a symmetric hyperbolic provided $\rho>0, \theta>0$.

\section{Local  existence}
In the following, we shall assume  that for $\theta>0$
\begin{align}\label{3.1}
a(\theta)>0, a^\prime(\theta)\ge0,  \frac{1}{2} \left(\frac{Z(\theta)}{\theta}\right)^\prime\ge0
\end{align}
The assumption $a^\prime(\theta)\ge0$ implies $e_\theta\ge C_v>0$, which make the system \eqref{1.1}-\eqref{1.3} uniformly hyperbolic without small condition. The third inequality in \eqref{3.1} will give the $L^2$  estimates of $q$  from Lemma \ref{le3.2} below, which will be used in the blow-up result.
Note also that by choosing $Z(\theta)=\frac{\tau_1(\theta)}{\kappa (\theta)}=k\theta^\alpha$ with $k$ be any constant and $1\le \alpha<2$, the assumption \eqref{3.1} holds.

Now, we transform the equations \eqref{1.1}-\eqref{1.3} into a first-order symmetric hyperbolic system. First, we rewrite the equation $\eqref{1.1}_3$  for $\theta$ as
\begin{align}
\rho e_\theta \theta_t+\left(\rho u e_\theta-\frac{2a(\theta)}{Z(\theta)}q\right)\theta_x+R \rho \theta u_x+q_x=\frac{2a(\theta)}{\tau_1(\theta)}q^2+\frac{1}{\mu}S^2.
\end{align}
Then, we have
\begin{align}
A^0(U) U_t+A^1(U) U_x+B(U)U=F(U),
\end{align}
where $U=(\rho, u, \theta, q, S)$ and
\begin{align*}
A^0(U)=\mathrm{diag} \{ \frac{R\theta}{\rho}, \rho, \frac{\rho e_\theta}{\theta} , \frac{\tau_1(\theta) \rho}{\kappa(\theta)}, \frac{\tau_2 \rho}{\mu}\},\\
A^1(U)=
\begin{pmatrix}
\frac{R\theta}{\rho} u& R \theta&0&0&0\\
R\theta& \rho u& R \rho& 0&-1\\
0& R\rho & \left( \frac{\rho u e_\theta}{\theta}-\frac{2a(\theta)}{\theta Z(\theta)}q\right)& \frac{1}{\theta}&0\\
0&0&\frac{1}{\theta}&\frac{\tau_1(\theta)}{\kappa(\theta)} \rho u &0\\
0&-1&0&0&\frac{\tau_2}{\mu} \rho u
\end{pmatrix},\\
B(U)=\mathrm{diag}\{0,0,0,\frac{1}{\kappa \theta}, \frac{1}{\mu}\}, F(U)=\mathrm{diag} \{ 0,0,-\frac{2 a(\theta)}{\tau_1(\theta)  \theta}q^2-\frac{S^2}{\mu \theta},0,0\}.
\end{align*}
Therefore, the local existence follows immediately, see  \cite{Kaw83,TAY,JR2000}.
\begin{theorem} \label{th2.1}
Let $s\ge 2$. Suppose that $$(\rho_0-1,u_0,\theta_0-1,q_0,S_0)\in H^s(\mathbb R)$$ with $\min_x (\rho_0(x), \theta_0(x))>0$,
there exists a unique local solution $(\rho,u,\theta,q,S)$ to \eqref{1.1}-\eqref{1.5} in some time interval $[0,T]$ with
\begin{align}
(\rho-1,u,\theta-1,q,S) \in C^0([0,T],H^s(\mathbb R)) \cap C^1([0,T],H^{s-1}(\mathbb R)),\\
\min_x (\rho(t,x), \theta (t,x))>0, \qquad \forall t>0.
\end{align}
\end{theorem}

\section{ Entropy  equation and global existence}
In this part, we first derive an entropy equation for system \eqref{1.1}-\eqref{1.3}.  Defining the entropy
\begin{align} \label{entropy}
\eta:=C_v  \ln \theta-R \ln \rho- \left(\frac{Z(\theta)}{2\theta}\right)^\prime q^2.
\end{align}
Similar to \cite{HRW}, we have for a local solution
\begin{lemma} \label{le3.1}
The entropy $\eta$ defined above satisfies
\begin{align}\label{2.1}
(\rho \eta)_t + \left(\rho u \eta+\frac{q}{\theta}\right)_x=\frac{q^2}{\kappa(\theta) \theta^2}+\frac{S^2}{\mu\theta}.
\end{align}
\end{lemma}
\begin{proof}
From the energy equation $\eqref{1.3}$, we easily get the equation for $e$ as follows:
\begin{align}
\rho e_t+\rho u e_x+p u_x+q_x=\frac{1}{\mu} S^2.
\end{align}
Dividing the above equation by $\theta$ and using formula \eqref{1.4}, one has
\begin{align*}
\frac{\rho}{\theta}(C_v \theta+a(\theta)q^2)_t+\frac{\rho u}{\theta} (C_v\theta+a(\theta) q^2)_x+R \rho u_x+\frac{q_x}{\theta}=\frac{1}{\mu \theta}S^2.
\end{align*}
Now, we calculate the following term
\begin{align*}
&\frac{\rho}{\theta}(a(\theta)q^2)_t+\frac{\rho u}{\theta} (a(\theta)q^2)_x\\
&=\rho\left(\frac{a(\theta)}{\theta} q^2\right)_t+\rho \frac{a(\theta)}{\theta^2}q^2 \theta_t+\rho u \left( \frac{a(\theta)}{\theta^2}q^2\right)_x +\rho u \frac{a(\theta)}{\theta} q^2 \theta_x\\
&=\rho\left(\frac{a(\theta)}{\theta} q^2\right)_t-\rho\left( \frac{Z(\theta)}{2\theta^2} q^2\right)_t+\rho \frac{Z(\theta)}{ \theta^2} q q_t+\rho u \left( \frac{a(\theta)}{\theta}q^2\right)_x-\rho u\left( \frac{Z(\theta)}{2\theta^2} q^2\right)_x+\rho u \frac{Z(\theta)}{\theta^2} q q_x\\
&=\rho \left( \left(\frac{a(\theta)}{\theta} -\frac{Z(\theta)}{2\theta^2}\right)q^2\right)_t+\rho u \left( \left(\frac{a(\theta)}{\theta} -\frac{Z(\theta)}{2\theta^2}\right)q^2\right)_x-\frac{1}{\kappa(\theta) \theta^2} q^2-\frac{\theta_x}{\theta^2} q
\end{align*}
where we have used the identity $(\frac{Z(\theta)}{2\theta^2})_t=-\frac{a(\theta)}{\theta^2} \theta_t$.
On the other hand, using the mass equation $\eqref{1.1}_1$, we have
\begin{align*}
R \rho u_x=-R \rho( (\ln \rho)_t+u(\ln \rho)_x)
\end{align*}
Combining the above estimates and noting that $$\frac{a(\theta)}{\theta}-\frac{Z(\theta)}{2\theta^2}=- \left(\frac{Z(\theta)}{2\theta}\right)^\prime,$$ we get the desired result.
\end{proof}
\begin{remark}
Once we use the constitutive relation \eqref{1.2}, there are three unknown thermodynamic variable (for example, the density, temperature and heat flux) rather than two in the  classical system. Thus, with the entropy defined in \eqref{entropy}, we can get an extended Gibbs relation used in  extended irreversible thermodynamics (see \cite{Rug83, Mu67}) as
\begin{align}\label{Gibb}
\theta \mathrm{d} \eta=\mathrm{d} e+p \mathrm{d} v-\frac{Z(\theta)}{\theta} q \mathrm{d} q,
\end{align}
where $v=\frac{1}{\rho}$ and $\eta$ is the physical entropy. When $\tau_1(\theta)=0$, the equation \eqref{Gibb}  reduces to the classical Gibbs relation.
\end{remark}
The entropy equation implies the following lower energy estimates:
\begin{lemma}\label{le3.2}
Let $(\rho, u, \theta, q, S)$ be local solutions to \eqref{1.1}-\eqref{1.5}, then we have
\begin{align}
\int_{\mathbb R} \left[ C_v \rho (\theta-\ln \theta-1)+R (\rho \ln \rho-\rho+1)+\rho a(\theta) q^2+\frac{\tau_2}{2\mu} S^2 +\frac{1}{2}\rho u^2\right]\dif x \nonumber\\
+\int_0^t \int_{\mathbb R} \left( \frac{q^2}{\kappa (\theta) \theta^2}+\frac{S^2}{\mu\theta}\right) \dif x \dif t =I_0, \label{3.6}
\end{align}
where
$$
I_0:= \int_{\mathbb R} \left(C_v \rho_0 (\theta_0-\ln \theta_0-1)+R(\rho_0 \ln \rho_0-\rho_0+1)+\rho_0 a(\theta_0) q_0^2+\frac{\tau_2}{\mu} S_0^2+\frac{1}{2}\rho_0 u_0^2 \right)\dif x.
$$
Moreover, for $|\rho-1|\le \frac{1}{2}, |\theta-1|\le \frac{1}{2}$, we have
\begin{align}\label{3.7}
\int_{\mathbb R} ((\theta-1)^2+(\rho-1)^2+q^2+S^2+u^2)\dif x+\int_0^t \int_{\mathbb R} (q^2+S^2)\dif x \dif t \le C I_0.
\end{align}
\end{lemma}
\begin{proof}
Combing the entropy equation \eqref{2.1}, the momentum equation $\eqref{1.1}_2$ and the energy equation $\eqref{1.1}_3$, we have
\small{
$$
 \left[C_v \rho (\theta-\ln \theta-1)+R(  \rho\ln \rho-\rho+1)+ \rho \left( a(\theta)+\frac{1}{2} \left(\frac{Z(\theta)}{\theta}\right)^\prime \right) q^2 +\frac{\tau_2}{2\mu}\rho S^2+\frac{1}{2} \rho u^2\right]_t +
 $$
 $$
 \left[C_v\rho u  (\theta-\ln \theta-1)+R \rho u \ln \rho -R \rho u+\rho u \left( a(\theta)+\frac{1}{2} \left(\frac{Z(\theta)}{\theta}\right)^\prime \right) q^2+\frac{\tau_2}{2\mu} \rho u S^2+\frac{1}{2} \rho u^3+pu+q-Su-\frac{q}{\theta}\right]_x
 $$
 $$
+\frac{q^2}{\kappa(\theta) \theta^2}+\frac{S^2}{\mu\theta}=0.  \nonumber
$$
}
Then, using assumption \eqref{3.1}, we get \eqref{3.6} immediately. Moreover, using Taylor expansions, we get
\begin{align}
\theta-\ln \theta-1=\frac{1}{2\xi^2} (\theta-1)^2, \\
\rho \ln \rho -\rho+1= \frac{1}{2 \eta} (\rho-1)^2,
\end{align}
where $ \xi\in (1,\theta), \eta \in (1,\rho)$. Therefore, by assuming $|\rho-1|\le \frac{1}{2}, |\theta-1|\le \frac{1}{2}$, we get  the $L^2$ estimates \eqref{3.7}.
\end{proof}
We give a remark on global existence for small data.
\begin{remark}
Since the system is symmetric hyperbolic, zero-order estimates (from the entropy equation) together with Kawashima's dissipation structure (from our linear analysis in \cite{HRW}: The linearized system are the same) would imply the global existence of strong solutions for small initial data which we do not state in detail.
\end{remark}

\section{Blow-up for large data}
 Here we now show that there exist  large initial data $(\rho_0, u_0, \theta_0, q_0, S_0)$ such that the local solution $(\rho, u ,\theta, q, S)(t,x)$ must blow up in finite time.

Since the system is symmetric hyperbolic, the local solutions of \eqref{1.1}-\eqref{1.3} possess  the finite propagation speed property:
\begin{proposition}\label{pro4.1}
Let $(\rho_0, u_0, \theta_0, q_0, S_0)$ be given as in Theorem \ref{th2.1} and $(\rho, u, \theta, q, S)$ be local solutions to \eqref{1.1}-\eqref{1.5} on  $[0, T_0)$. Let $M>0$. We further assume the initial data
$(\rho_0-1, u_0, \theta_0-1, q_0, S_0)$ are compactly supported in $(-M, M)$.
Then, there exists a constant $\sigma$ such that
\begin{align*}
(\rho(\cdot,t), u(\cdot, t), \theta(\cdot,t), q(\cdot, t), S(\cdot, t)=(1, 0, 1, 0, 0):=(\bar \rho, \bar u ,\bar \theta, \bar q, \bar S)
\end{align*}
on $D(t)=\{ x \in \mathbb R \big{|} |x|\ge M+\sigma t\},\quad 0\le t < T_0.$
\end{proposition}
Now, we define some averaged quantities as follows:
\begin{align}
F(t):=\int \rho u \cdot x \dif x -\tau_2 \int \rho S \dif x,\\
G(t):=\int_{\mathbb R} (E(x,t)-\bar E) \dif x,
\end{align}
where
$$
 E=\frac{1}{2} \rho u^2 +\frac{\tau_2}{2\mu} \rho S^2+\rho e(\theta, q)
$$
is the total energy and
$$
\bar E:=\bar \rho(\bar e+\frac{1}{2} \bar u ^2)=C_v.
$$
 The functional $F$ with the second term involving $S$ is different from those used in \cite{S85,HRW}.

We mention that the functional defined above exists since the solution $(\rho-1, u, \theta-1, q, S)$ is zero on the set $D(t)$ defined in Proposition \ref{pro4.1}.

Now, we are ready to show our main result.
\begin{theorem}\label{th1}
We assume that the initial data satisfy the assumption in  Theorem \ref{th2.1} and Proposition \ref{pro4.1} . Moreover, we assume that  assumption \eqref{3.1} holds and
\begin{equation} \label{assongnull}
G(0)>0.
\end{equation}
Then, there exists $(\rho_0, u_0, \theta_0, q_0, S_0)$ satisfying
\begin{align}
F(0)>\frac{32 \sigma \max \rho_0}{3-\gamma} M^2
\end{align}
and
\begin{align}
4\left(\frac{(3-\gamma) \mu \tau_2}{M^2}+\gamma-1\right) (H_0+\frac{\max \rho_0}{2} \|u_0\|_{L^2}^2 ) \le \frac{128 \sigma^2 \max \rho_0 M}{3-\gamma},
\end{align}
where $H_0$ is defined in \eqref{H0},
such that the length $T_0$ of the maximal interval of existence of a smooth solution $(\rho, u, \theta, q, S)$ of \eqref{1.1}-\eqref{1.5} is finite,
provided the compact support of the initial data is sufficiently large and $\gamma:= 1 + \frac{R}{C_v} $is sufficiently close to $1$.
\end{theorem}
\begin{proof}
From equations $\eqref{1.1}_2$ and $\eqref{1.1}_3$, we can get the equation for $E$:
\begin{align*}
E_t+(uE+up-uS+q)_x=0,
\end{align*}
which implies that $G$ is a constant and
\begin{align}
G(t)=G(0)>0.
\end{align}
In the following, $\int$ denotes $\int_{\mathbb R}$ for simplicity.
Using the momentum equation $\eqref{1.1}_2$,  the constitutive equation \eqref{1.3}, Lemma \ref{le3.2} and \eqref{3.1}, we can derive
\begin{align*}
F^\prime(t)&=\int (\rho u)_t \cdot x \dif x -\tau_2 \int (\rho S)_t \dif x\\
&=\int (-\rho u^2 -p+S)_x\cdot x \dif x +\int S \dif x\\
&=\int \rho u^2 +\int (p-\bar p)\dif x\\
&=\int \rho u^2 \dif x+\int (R\rho \theta-R \bar \rho \bar \theta) \dif x\\
&=\int \rho u^2 \dif x+ \int \left( \frac{R}{C_v} (\rho e-\bar \rho \bar e)- \frac{R}{C_v} a(\theta) \rho q^2-\frac{R}{C_v} \frac{\tau_2}{2\mu}\rho S^2 \right)\dif x\\
&=\int \rho u^2 \dif x+(\gamma-1) \int (E-\bar E) \dif x - (\gamma-1)\int \frac{1}{2} \rho u^2 \dif x -(\gamma-1) \int \left(a(\theta) \rho q^2+\frac{\tau_2}{2\mu} \rho S^2\right) \dif x\\
&\ge \frac{3-\gamma}{2} \int \rho u^2 \dif x-(\gamma-1) \int  \left( a(\theta)\rho q^2+\frac{\tau_2}{2\mu} \rho S^2\right) \dif x\\
&\ge \frac{3-\gamma}{2} \int \rho  u^2 \dif x -(\gamma-1) (H_0+\frac{\max \rho_0}{2} \|u_0\|_{L^2}^2 )
\end{align*}
where  $\gamma=\frac{R}{C_v}+1$ and
\begin{align}\label{H0}
H_0:=\int C_v \rho_0(\theta_0-\ln \theta_0-1)+R(\rho_0 \ln \rho_0-\rho_0+1)+\rho_0\left(a(\theta_0)+\frac{1}{2} \left(\frac{Z(\theta_0)}{\theta_0}\right)^\prime\right) q_0^2+\frac{\tau_2}{2\mu} S_0^2 \dif x.
\end{align}

On the other hand,
\begin{align*}
F^2(t) &\le 2\left( \int \rho u \cdot x \dif x\right)^2+2\tau_2^2 \left(\int \rho S \dif x\right)^2\\
&\le 4 \max \rho_0 (M+\sigma t)^3 \int \rho u^2 \dif x+2\tau_2^2 \int \rho S^2 \dif x \int \rho \dif x\\
&\le 4\max \rho_0 (M+\sigma t)^3 \int \rho u^2 \dif x+4\mu \tau_2 (H_0+\frac{\max \rho_0}{2} \|u_0\|_{L^2}^2) \max \rho_0 \cdot 2 (M+\sigma t)
\end{align*}
which implies
\begin{align*}
\int \rho u^2 \dif x \ge \frac{F(t)^2}{4\max \rho_0 (M+\sigma t)^3} -\frac{2 \mu \tau_2 (H_0+\frac{\max \rho_0}{2} \|u_0\|_{L^2}^2)}{(M+\sigma t)^2}.
\end{align*}
Combining the above results, we derive
\begin{align}
F^\prime(t) &\ge \frac{3-\gamma}{8 \max \rho_0 (M+\sigma t)^3}F^2(t)-\left(\frac{(3-\gamma)\mu \tau_2}{ (M+\sigma t)^2}+\gamma-1\right)(H_0+\frac{\max \rho_0}{2} \|u_0\|_{L^2}^2)\\\nonumber
&\equiv \frac{c_3}{(1+c_2 t)^3}F(t)^2-K(t)
\end{align}
where $c_2:= \frac{\sigma}{M}, c_3:=\frac{3-\gamma}{8\max \rho_0 M^3}$.
With this Riccati inequality, we can  show the blow-up result.

Indeed, assuming a priori that
\begin{align}\label{AP1}
2 K(t)\le\frac{c_3}{(1+c_2 t)^3}F^2(t),
\end{align}
then we have
\begin{align*}
F^\prime (t) \ge \frac{c_3}{2(1+c_2 t)^3} F^2(t),
\end{align*}
which gives
\begin{align}
\frac{1}{F_0}\ge \frac{1}{F_0}-\frac{1}{F(t)} \ge \frac{c_3}{4c_2}-\frac{c_3}{4c_2 (1+c_2 t)^2}
\end{align}
for which the maximal existence time $T$ can not be infinity provided
\begin{align}\label{AS1}
F_0> \frac{4c_2}{c_3} =\frac{32\sigma \max \rho_0 M^2}{3-\gamma}.
\end{align}
Here $F_0=F(0)$.
Moreover, we have
\begin{align}
\frac{1}{F(t)} \le \frac{1}{F_0}-\frac{c_3}{4c_2} +\frac{c_3}{4c_2 (1+c_2t)^2},
\end{align}
which implies that
\begin{align}\label{4.12}
F(t)\ge \frac{4 c_2 (1+c_2 t)^2}{c_3}.
\end{align}
 To show that the a priori estimate \eqref{AP1} holds, we use the bootstrap method expressed in the following simple lemma.
\begin{lemma} \label{le4.3}
Let $f \in C^0\left([0,\infty),[0,\infty)\right)$ and $0<a<b$ such that the following holds for any $0\leq \alpha < \beta < \infty$:
$$
f(0) < a \qquad  \mbox{\rm and} \qquad
\left(
\forall \, t \in [\alpha,\beta]: f(t) \leq b\; \implies \; \forall \, t \in [\alpha,\beta]: f(t) \leq a.
 \right).
 $$
Then we have
$$
\forall t\geq 0:\; f(t) \leq a.
$$
\end{lemma}
 That is, under the a priori assumption \eqref{AP1}, we need to show that
\begin{align}\label{AP2}
4 K(t) \le \frac{c_3}{(1+c_2 t)^3}F^2(t).
\end{align}
We need the above equality holds  in particular for $t=0$, that is,
\begin{align}\label{AS2}
4 \left(\frac{(3-\gamma) \mu \tau_2}{M^2}+\gamma-1\right) (H_0+\frac{\max \rho_0}{2} \|u_0\|_{L^2}^2 ) \le c_3 F_0^2.
\end{align}
Next, using \eqref{4.12}, one only need to show
\begin{align}
4 K(t) \frac{(1+c_2 t)^2}{c_3} \le \frac{16 c_2^2}{c_3^2} (1+c_2t)^4
\end{align}
which is sufficient to show
\begin{align}\label{AS3}
4\left(\frac{(3-\gamma) \mu \tau_2}{M^2}+\gamma-1\right) (H_0+\frac{\max \rho_0}{2} \|u_0\|_{L^2}^2 ) \le \frac{16 c_2^2}{c_3}
\end{align}
Note that \eqref{AS1} and \eqref{AS3} imply \eqref{AS2}, we need to find some $u_0$ such that the assumptions \eqref{AS1} and \eqref{AS3} hold.
As in \cite{HR2}, we choose $u_0\in H^2(\mathbb R)\cap  C^1(\mathbb R)$ as follows:
\begin{equation}
u_0(x):=
\begin{cases}
0,\qquad \qquad\qquad \qquad\qquad \quad x\in(-\infty,-M],\\
\frac{L}{2}\cos(\pi (x+M))-\frac{L}{2}\qquad \,x\in(-M,-M+1],\\
-L,\qquad \qquad\qquad\qquad \qquad x\in(-M+1,-1],\\
L\cos(\frac{\pi}{2}(x-1)),\qquad \qquad\, x\in(-1,1],\\
L,\qquad \qquad\qquad  \qquad \qquad x\in (1,M-1],\\
\frac{L}{2}\cos(\pi(x-M+1))+\frac{L}{2}\quad x\in(M-1,M],\\
0,\qquad\qquad \qquad\qquad \qquad  x\in(M,\infty),
\end{cases}
\end{equation}
where $L$ is a positive constant to be determined later. We assume
$M\geq 4$.
 Assumption (\ref{assongnull}) can easily be satisfied since it is equivalent to requiring
$$
\int_{\mathbb R}\left(\rho_0e_0 -\bar \rho \bar e + \frac{1}{2} u_0^2\right) dx > 0,
$$
which is satisfied by choosing $\rho_0\theta_0 >  \bar \rho \bar \theta=1$.
Since
\begin{align*}
\int_{\mathbb R}(x\rho_0(x) u_0(x))\dif x
 \ge \frac{L}{2}\min\rho_0M^2
 \end{align*}
 and
 \begin{align*}
 \left| \tau_2 \int \rho_0 S_0 \dif x\right|\le \int_{-M}^M \rho_0 \dif x+\tau_2 \int \rho_0 S_0^2 \dif x\le \max \rho_0 (1+\mu H_0^2) M^2.
 \end{align*}
We choose $L$ large
 enough, and independent of $M$,
  such that
  $$\frac{L}{4}\mathrm{min} \rho_0 >\max \{ \max \rho_0(1+\mu H_0^2), \frac{32\sigma \max \rho_0}{3-\gamma}\}$$
  Therefore,  \eqref{AS1} hold. Now,  after having chosen $\sigma$ large enough,  fix $L$. Then we choose $M$ sufficiently large and $\gamma-1$ sufficiently small such that \eqref{AS3} holds.
\end{proof}
{\bf Acknowledgement:} Yuxi Hu's Research is supported by the Fundamental Research Funds for the Central Universities (No. 2023ZKPYLX01).

\end{document}